\documentclass[12pt, a4paper]{article}
\usepackage[latin1]{inputenc}
\usepackage[english]{babel}
\usepackage{indentfirst}
\usepackage{amstext}
\usepackage{setspace}
\usepackage{amsfonts}
\usepackage{textcomp}
\usepackage{amssymb}
\usepackage{amscd}
\usepackage{epsf}
\usepackage{graphicx}
\usepackage{epsfig}

\newcommand {\demo}{\hskip -0.6cm{\bf Proof:  }}
\newcommand {\fim}{\hfill{$\square$}\vskip 1pc}

\newcommand {\N}{\mathbb{N}}
\newcommand {\Z}{\mathbb{Z}}

\newcommand {\F}{\mathbb{F}}

\newtheorem{theorem}{Theorem}[section]
\newtheorem{lema}[theorem]{Lemma}
\newtheorem{corolario}[theorem]{Corollary}
\newtheorem{definition}[theorem]{Definition}
\newtheorem{proposition}[theorem]{Proposition}

\newtheorem{obs}[theorem]{Remark}

\begin{document}

\title{Leavitt path algebras as partial skew group rings}
\maketitle
\begin{center}
%{\Large{\bf Perron-Frobenius operators and representations of the Cuntz-Krieger algebras for infinite matrices}} \\
%\vspace{5mm}
{\large Daniel Gonçalves and Danilo Royer}\\
\end{center}  
\vspace{8mm}

\abstract 
We realize Leavitt path algebras as partial skew group rings and give new proofs, based on partial skew group ring theory, of the Cuntz-Krieger uniqueness theorem and simplicity criteria for Leavitt path algebras. 
\doublespace

\section{Introduction}

Leavitt path algebras were introduced by Abrams and Aranda Pino as a generalization of Leavitt algebras and as algebraic analogues of Cuntz-Krieger C*-algebras (see \cite{arandapino,arandapino1,CK}). Ever since their introduction, Leavitt path algebras kept growing in importance, as researchers characterized their properties (such as simplicity, finite dimensionality, locally finiteness, etc.) in terms of intrinsic properties of the underlying graph, (see \cite{arandapinomolina, AALP, AraPardo, Hazrat}), and studied their relations with other fields of mathematics, as for example operator algebras theory (see \cite{tomforde}).

It is our intention to connect the theory of Leavitt path algebras with another key concept arising from operator algebras theory, namely, partial skew group rings (or just skew group rings, depending on the literature, see \cite{Avilaferrero, Ex, Ferrerolazzarin}). 

Partial skew group rings (as we will call it) were introduced by Exel and Dokuchaev in \cite{Ex} as algebraic analogues of C* partial crossed products. The later algebras proved to be a powerful tool in the study of operator algebras, (see \cite{Ex1, Ex2, Ex3, QR}) and so it is important to realize C*-algebras as partial crossed products (see \cite{BE, EL} for example), as one can then benefit from the established theory about partial crossed products. 

The goal of our paper is to allow the theory of Leavitt path algebras to benefit from the theory of partial skew group rings and vice versa. For this we realize Leavitt path algebras as partial skew group rings (section 3) and, to illustrate how this interaction can be fruitful,  in section 4 we provide new proofs of the Cuntz-Krieger uniqueness theorem and the simplicity criteria for Leavitt Path Algebras (originally proved in \cite{arandapino, tomforde}). We remark that the proofs we provide for these results rely only on the theory of partial skew group rings. Furthermore, from the partial skew group ring point of view, the simplicity criteria for Leavitt path algebras can be seen as an analogue for the simplicity criteria for certain partial crossed products described in \cite{ELQ} (but the results in \cite{ELQ} are more general, so, to our knowledge, a full algebraic analogue has not been discovered yet). 

To define the relevant partial skew group ring we need first to define the appropriate partial action. We do this in section 2 below.

%it is interesting to notice that much as it happened with Leavitt path algebras, some of the results for C*-algebraic partial crossed products have been generalized to the algebraic setting and other independent results have been found (see \cite{Ex, Avilaferrero, Ferrerolazzarin}). 
%In this direction, 

\section{The partial action}
\label{section2}

In this section we associate to a graph a partial action of the free group on the edges of the graph. For the definition of partial actions and partial skew group rings we refer the reader to \cite{Ex}.

%reader convenience we recall the definitions of a partial actions of algebras and partial skew group rings below.
%\begin{definition} Let $G$ be a group and $A$ an algebra. A partial action $\alpha$ of $G$ on $A$ is a collection of ideals $D_g$, $g\in G$, and algebra isomorphisms $\alpha_g : D_{g^-1} \rightarrow D_g$ such that:
%\begin{itemize}
%\item  $D_1 = A$ and $\alpha_1$ is the identity map in $A$.
%\end{itemize}
%\end{definition}

Let $E=(E^0, E^1, r, s)$ be a directed graph. We denote the set of all finite paths in $E$ by $W$, and the set of all infinite paths in $E$ by $W^\infty$. In symbols we have:
$$W:=\bigcup\limits_{n=1}^\infty \left\{\xi_1...\xi_n:\,\, \xi_i\in E^1 \text{ and } r(\xi_i)=s(\xi_{i+1}) \forall i\in\{1,...,n-1\}\right\} \text{ and } $$
$$W^\infty:=\left\{\xi_1\xi_2\xi_3... : \,\,\xi_i\in E^1\,\, \forall i\in \N \text{ and } r(\xi_i)=s(\xi_{i+1})\forall i\in \N\right\}.$$
Notice that $W$ ($W^\infty$) is a subset of the set of all finite (infinite) words in the alphabet $E^1$. 

As usual, we extend the range and source maps to $W\cup W^\infty\cup E^0$ by defining $s(\xi):=s(\xi_1)$, for $\xi=\xi_1 \xi_2... \in W^\infty$ or $\xi=\xi_1\ldots \xi_n \in W$,  $r(\xi):=r(\xi_n)$ for $\xi=\xi_1\ldots \xi_n \in W$ and $r(v)=s(v)=v$ for $v\in E^0$. We say that a path $\xi = \xi_1 \ldots \xi_n$ is the beginning of the path $\eta=\eta_1 \ldots \eta_m$ if $m \geq n$ and $\xi_i = \eta_i$ for $i=1,..,n$.

Our goal is to define a partial action associated to the free group, $\F$, generated by $E^1$. The first obstacle we have to overcome is to decide over which set $\F$ should act. Studying similar problems in C*-algebra theory (see for example \cite{GRHouston}) we notice that the best candidate of this set is the spectrum of the algebra, and this leads us to define:
$$X=\left\{\xi\in W : r(\xi) \text{ is a sink} \right\}\cup\{v\in E^0: v \text{ is a sink }\}\cup W^\infty.$$

The next task in our hands is to define the subsets $X_c$, $c\in \F$. For this, notice that we can see $W$ as a subset of $\F$, and define:

\begin{itemize}
\item $X_0:=X$, where $0$ is the neutral element of $\F$. 

\item $X_{b^{-1}}:=\{\xi\in X: s(\xi)=r(b)\},$ for all $b\in W$. 

\item $X_a:=\{\xi\in X: \xi_1\xi_2...\xi_{|a|}=a\},$ for all $a\in W$, where $|a|$ denotes the length of $a$, that is, if $a=a_1...a_n$, then $|a|=n$.

\item $X_{ab^{-1}}:=\{\xi\in X: \xi_1\xi_2...\xi_{|a|}=a\}=X_a,$ for $ab^{-1}\in \F$ with $a,b\in W$, $r(a)=r(b)$ and $ab^{-1}$ in its reduced form (that is, $a_{|a|}\neq b_{|b|}$).

\item $X_c:=\emptyset$, for all other $c \in \F$.

\end{itemize}

\begin{obs} In the definitions above, note that $r(b)\in X_{b^{-1}}$ if and only if $r(b)$ is a sink, and moreover, if $r(b)$ is a sink, then $X_{b^{-1}}=\{r(b)\}$ and $X_b=\{b\}$.
\end{obs}

%Above we have defined all subsets of $X$ associated to elements of $\F$. Next we will describe some of the relations between these sets and their relations with sets of the form $X_v:=\{\xi\in X:s(\xi)=v\}$ , $v\in E^0$ (Notice that since $s(v)=v$, then $v\in X_v$ if and only if $v$ is a sink.). 
For future reference we describe below some of the relations between the sets $X_c$, $c\in\F$, defined above and their relations with sets of the form $X_v=\{\xi\in X:s(\xi)=v\}$, $v\in E^0$. (Notice that since $s(v)=v$, then $v\in X_v$ if and only if $v$ is a sink, and in this case, $X_v=\{v\}$).

\begin{lema}\label{lema1}
Let $a,c\in W$, $b,d\in W\cup\{0\}$ and $v\in E^0$. Then: 

\begin{enumerate}

\item $X_{a^{-1}}\cap X_{c^{-1}}=\left\{\begin{array}{ll} 
X_{a^{-1}}=X_{c^{-1}}, & \text{ if } $r(a)=r(c)$\\
\emptyset , & \text{ otherwise }
\end{array}
\right ..$

\item $X_{a^{-1}}\cap X_{cd^{-1}}=\left\{\begin{array}{ll} 
X_{cd^{-1}}, & \text{ if } $r(a)=s(c)$\\
\emptyset ,& \text{ otherwise }
\end{array}
\right ..$

\item $X_{ab^{-1}}\cap X_{cd^{-1}}=\left\{\begin{array}{ll} 
X_{ab^{-1}}, & \text{ if } a=ct \text{ for some } t\in W\cup \{0\} \\
X_{cd^{-1}}, & \text{ if } c=at \text{ for some } t\in W \cup \{0\} \\
\emptyset , & \text{ otherwise }
\end{array}
\right .,$
 (supposing $r(a)=r(b)$ and $r(c)=r(d)$).
\item $X_v\cap X_{b^{-1}}=\left\{\begin{array}{ll} 
X_v=X_{b^{-1}}, & \text{ if } $r(b)=v$\\
\emptyset , & \text{ otherwise }
\end{array}
\right ..$

\item $X_v\cap X_{ab^{-1}}=\left\{\begin{array}{ll} 
X_{ab^{-1}}, & \text{ if } $s(a)=v$\\
\emptyset , & \text{ otherwise }
\end{array}
\right .$ and $X_v=\displaystyle \bigcup_{s(a)=v} X_{ab^{-1}}$.

\end{enumerate}
\end{lema}

\demo
These properties follow straightforward from the definitions of the sets.
\fim

The last thing we need to define, to obtain a partial action of $\F$ on $X$, are the bijective maps $\theta_c:X_{c^{-1}}\rightarrow X_c$, for $c\in \F$ such that $X_{c^{-1}}\neq \emptyset$, that is, for $c$ of the form $c=ab^{-1}$, with $a,b\in W\cup\{0\}$.

As required by the notion of a partial action we define $\theta_{0}:X_{0}\rightarrow X_{0}$ as the identity map. For $b\in W$, we define $\theta_b:X_{b^{-1}}\rightarrow X_b$ as the ``add b" or ``creation" map, that is, if $\xi\in X_{b^{-1}}$ then $\theta_b(\xi)=b\xi$. Notice that $\theta_b(\xi)$ is well defined, since $\xi\in X_{b^{-1}}$ implies that $s(\xi)=r(b)$ and so  $b\xi\in X_b$. Also, notice that in case $r(b)$ is a sink then $X_{b^{-1}}=\{r(b)\}$, and we are assuming above that $br(b)=b$.

 It is not hard to see that $\theta_b$ is bijective, with inverse $\theta_{b^{-1}}:X_b\rightarrow X_{b^{-1}}$ given by the ``erase b" map, that is,
$\theta_{b^{-1}}(\eta)= \eta_{|b|+1}\eta_{|b|+2}...$ if $r(b)$ is not a sink  and $\theta_{b^{-1}}(b)=r(b)$, if $b$ is a sink. 
 
Now, for $a,b\in W$ with $r(a)=r(b)$ and $a_{|a|}\neq b_{|b|}$ we define $\theta_{ab^{-1}}:X_{ba^{-1}}\rightarrow X_{ab^{-1}}$ as the ``erase b and add a" map, that is, $\theta_{ab^{-1}}(\xi)=a\xi_{(|b|+1)}\xi_{(|b|+2)}...$. Notice that $\theta_{ab^{-1}}(\xi)$ is well defined, since if $\xi\in X_{ba^{-1}}$ then $\xi_1...\xi_{|b|}=b$ and since $r(a)=r(b)=s(\xi_{|b|+1})$ we have that $a\xi_{|b|+1}\xi_{|b|+2}...\in X_{ab^{-1}}$.

Again, it is not hard to see that $\theta_{ab^{-1}}$ is bijective with inverse $\theta_{ba^{-1}}:X_{ab^{-1}}\rightarrow X_{ba^{-1}}$ given by $\theta_{ba^{-1}} (\eta)=b\eta_{(|a|+1)}\eta_{(|a|+2)}...$ .

We have now defined all the necessary ingredients for a partial action of sets in $X$. We leave to the reader the standard verification that $\{\{X_c\}_{c\in \F}, \{\theta_c\}_{c\in \F}\}$ satisfy the axioms of a partial action of sets as defined in \cite{Ex}.

The reader might be wandering how the above partial action relate with our work, since our aim is to realize the Leavitt path algebra of a graph as a partial skew group ring and so far we only have a partial action on the set level. In the remainder of this section we show how the set partial action defined above induces a partial action in the algebra level (for more details about the relations between partial actions on sets and partial actions of algebras see \cite{Vivi}).

 Let $F(X)$ be the $K$-algebra of the functions from $X$ to $K$ with pointwise multiplication. For each $c\in \F$, with $X_c\neq\emptyset$, let $F(X_c)$ be the $K$-algebra of functions from $X_c$ to $K$. Note that $F(X_c)$ may be identified with the subset of the functions in $F(X)$ that vanishes outside of $X_c$. For $c\in \F$, with $X_c=\emptyset$, let $F(X_c)$ be the subset of $F(X)$ that contains only the null function. Note that each $F(X_c)$ is an ideal of the $K$-algebra $F(X)$.

Now, for each $c\in\F$, define $\alpha_c:F(X_{c^{-1}})\rightarrow F(X_c) $ by $\alpha(f) = f\circ \theta_{c^{-1}}$, which is an $K$-isomorphism. One can now check that the family $\{\{\alpha_c\}_{c\in \F}, \{F(X_c)\}_{c\in \F}\}$ is a partial action of $\F$ on $F(X)$ (analogous verifications are done in \cite{Vivi} for a similar case).

Even though we have obtained a partial action in the algebra level, the associated partial skew group ring is too ``large" for our purposes. To obtain the Leavitt path algebra we need to make our algebra ``smaller". We do this by choosing appropriate functions in $F(X)$, that is, by letting the ideals in the definition of the partial action to be the linear span of some characteristic functions. We make this precise below.

For each $c\in \F$, let $1_c$ be the characteristic function of the set $X_c$, that is, $1_c(\xi)=1_{X_c}(\xi)=[\xi\in X_c]$, where $[\xi\in X_c]=1$, if $\xi\in X_c$ and $[\xi\in X_c]=0$, if $\xi\notin X_c$ and let $1_v$ denote the characteristic function of the set $X_v$.

 \begin{lema}\label{lema2} Let $p,q\in \F$. Then 
 \begin{enumerate}
 \item $\alpha_p(1_{p^{-1}}1_q)=1_p1_{pq}$. 
 \item For $a\in W$ and $b\in W\cup\{0\}$, it holds that $$\alpha_a(1_{a^{-1}}1_v)=\left\{\begin{array}{ll} 1_a, & \text{ if } r(a)=v\\
 0, & \text{ otherwise }
\end{array} \right .$$
 and $$\alpha_{ab^{-1}}(1_{ba^{-1}}1_v)=\left\{\begin{array}{ll} 1_{ab^{-1}}, & \text{ if } s(b)=v\\
 0, & \text{ otherwise }
 \end{array}\right .
$$
 \end{enumerate}
 \end{lema}

 \demo 
This lemma is essentially a consequence of the definitions of $X_p$ and lemma \ref{lema1}. We exemplify the techniques involved by giving the proof of the first statement when $p=ba^{-1}$ and $q=cd^{-1}$, with $a,c\in W$, $b,d\in W$ and with $r(a)=r(b)$ and $r(c)=r(d)$.
 
So, first suppose $c=at$, with $t\in W\cup\{0\}$. Then, for all $\xi \in X$,

 $\alpha_p(1_{p^{-1}}1_q)(\xi)= \alpha_{ba^{-1}}(1_{ab^{-1}}1_{cd^{-1}})(\xi)= \alpha_{ba^{-1}}(1_{ab^{-1}})(\xi)\alpha_{ba^{-1}}(1_{cd^{-1}})(\xi)$

$= 1_{ab^{-1}}(\theta_{a b^{-1}}(\xi)) 1_{atd^{-1}}(\theta_{a b^{-1}}(\xi)) = [\xi \in X_{bt}] = [\xi \in X_{btd^{-1}}] $

$= 1_{ba^{-1}}(\xi) 1_{btd^{-1}}(\xi) = 1_p(\xi) 1_{pq}(\xi) $, as desired.

Proceeding analogously we have that the first statement of the lemma is also valid if $a=ct$ for some $t\in W$. Finally, if $a\neq ct$, for all $t\in W$ and $c\neq at$, for all $t\in W$, then $1_{ab^{-1}}1_{cd^{-1}}=0$ by lemma \ref{lema1} and $1_{ba^{-1}cd^{-1}}=0$ by definition, and so the desired equality also holds. The other possibilities for $p, q \in \F$ are handled in a similar way and we ommit the proof in this cases.

As for the second statement of the lemma, this is a straightforward consequence of items 4. and 5. of lemma \ref{lema1}.

% By lemma \ref{lema1}, $X_{ab^{-1}}\cap X_{cd^{-1}}=X_{cd^{-1}}$, and so $1_{ab^{-1}}1_{cd^{-1}}=1_{cd^{-1}}$.
 
%Then, for each $\xi\in X_{ba^{-1}}$, 
%$$\alpha_{ba^{-1}}(1_{ab^{-1}}1_{cd^{-1}})(\xi)=\alpha_{ba^{-1}}(1_{cd^{-1}})(\xi)=1_{cd^{-1}}(\theta_{ab^{-1}}(\xi))=$$
%$$=[a\xi_{|b|+1}\xi_{|b|+2}...\in X_{cd^{-1}}]=[\xi_{|b|+1}\xi_{|b|+2}...\xi_{|b|+|t|}=t]=$$
%$$=[\xi_1\xi_2...\xi_{|b|+|t|}=bt]=[\xi\in X_{btd^{-1}}]=1_{btd^{-1}}(\xi)=1_{ba^{-1}cd^{-1}}(\xi).$$ If $\xi\notin X_{ba^{-1}}$ then $$\alpha_{ba^{-1}}(1_{ab^{-1}1_{cd^{-1}}})(\xi)=0=1_{ba^{-1}}(\xi)1_{ba^{-1}cd^{-1}}(\xi).$$ So, it follows that $\alpha_{ba^{-1}}(1_{ab^{-1}}1_{cd^{-1}})=1_{ba^{-1}}1_{ba^{-1}cd^{-1}}$. The same equality holds if $a=ct$ for some $t\in W$. If $a\neq ct$ for all $t\in W$ and $c\neq at$ for all $t\in W$ then $1_{ab^{-1}}1_{cd^{-1}}=0$ by lemma \ref{lema1} and $1_{ba^{-1}cd^{-1}}=0$ by definition, an so the desired equality also holds.
\fim

We can now define the partial action which induces the partial skew group ring that is isomorphic to the Leavitt path algebra.

Let $$D(X)=D_0=span\{\{1_p:p\in \F\setminus\{0\}\}\cup\{1_v:v\in E^0\}\},$$  (where {\it span} means the $K$-linear span) and, for each $p\in \F\setminus\{0\}$, let $D_p\subseteq F(X_p)$ be defined as $1_p D_0$, that is, $$D_p=span\{\{1_p1_q:q\in \F\}\}.$$

By lemma \ref{lema1}, $D(X)$ and $D_p$ are $K$-algebras, and moreover, $D_p$ is an ideal of $D$, for each $p\in \F$. Furthermore, since $\alpha_p(1_{p^{-1}}1_q)=1_p1_{pq}$, consider, for each $p\in \F$, the restriction of $\alpha_p$ to $D_{p^{-1}}$. Notice that $\alpha_{p}:D_{p^{-1}}\rightarrow D_p$ is an isomorphism of $K$-algebras. This way a partial action, in the sense of \cite{Ex}, is obtained, namely, $\{\{\alpha_p\}_{p\in \F}, \{D_p\}_{p\in \F}\}$. Let $D(X)\rtimes_\alpha\F$ be the partial skew group ring associated to it.

\begin{obs} Notice that, by lemma \ref{lema2}, the action $\alpha$ is independent of the action $\theta$ when restricted to the ideals $D_p$. So any set level partial action that implements an isomorphic correspondence between the sets $X_c$, $c\in \F$, and respects the correspondence between these sets given by $\theta$, will induce the same partial action $\{\{\alpha_p\}_{p\in \F}, \{D_p\}_{p\in \F}\}$ in the algebra level.
\end{obs}

\section{The Leavitt path algebra and skew group rings}

In this section we will show that, given a graph, the associated Leavitt path algebra is isomorphic to the partial skew group ring defined above.

For the reader´s convenience, we first recall the definition of Leavitt path algebras. More details can be found in \cite{arandapino, gonroy}.

\begin{definition} Let $E$ be a directed graph, and $K$ be a field. The Leavitt path algebra of $E$ with coefficients in $K$, denoted by $L_K(E)$, is the universal $K$-algebra generated by a set $\{v: v \in E^0\}$, of pairwise orthogonal idempotents, together with a set $\{e, e^* : e \in E^1\}$ of elements satisfying:
\begin{enumerate}
\item $s(e)e = er(e) = e$ for all $e\in E^1$,
\item $r(e)e^* = e^*s(e) = e^*$ for all $e\in E^1$,
\item for all $e,f\in E^1$, $e^*f =0$ if $e\neq f$ and $e^*e=r(e)$. 
\item $v =\sum\limits_{e\in E^1:s(e)=v}ee^*$ for every vertex $v$ with $0 < \#\{e: s(e) = v\} <\infty.$
\end{enumerate}
\end{definition}

\begin{proposition} There exists a K-homomorphism $\varphi$, from $L_K(E)$ on $D(X)\rtimes_\alpha\F,$ such that $\varphi(e)=1_e\delta_e$, $\varphi(e^*)=1_{e^{-1}}\delta_{e^{-1}}$, for all $e\in E^1$, and $\varphi(v)=1_v\delta_0$, for all $v\in E^0$.
\end{proposition}

\demo Consider the sets $\{1_e \delta_e, 1_{e^{-1}} \delta_{e^{-1}}: e\in E^1\}$ and $\{1_v\delta_0: v\in E^0\}$ in $ D(X)\rtimes_\alpha \F$. We will show that these sets satisfy the relations defining the Leavitt path algebra and then use the universal property of $L_K(E)$ to obtain the desired homomorphism.

First notice that, for each $e\in E^1$, $1_{s(e)}\delta_01_e\delta_e=1_{s(e)}1_e\delta_e$. By lemma \ref{lema1}, $1_{s(e)}1_e=1_e$, and so $1_{s(e)}1_e\delta_e=1_e\delta_e$. Also, $1_e\delta_e1_{r(e)}\delta_0=\alpha_e(\alpha_{e^{-1}}(1_e)1_{r(e)})\delta_e=\alpha_e(1_{e^{-1}}1_{r(e)})\delta_e.$
By lemma \ref{lema2}, $\alpha_e(1_{e^{-1}}1_{r(e)})=1_e$, and so $1_e\delta_e1_{r(e)}\delta_0=1_e \delta_e$ as desired. 

Relation two in the definition of the Leavitt algebra follows analogously to what is done above. Next we prove condition three. 

Let $f,g\in E^1$. Then $1_{e^{-1}}\delta_{e^{-1}}1_f\delta_f=\alpha_{e^{-1}}(\alpha_e(1_{e^{-1}})1_f)\delta_{e^{-1}f}=\alpha_{e^{-1}}(1_e1_f)\delta_{e^{-1}f}$. By lemma \ref{lema1}, if $e\neq f$ then $1_e 1_f=0$ and hence $1_{e^{-1}}\delta_{e^{-1}}1_f\delta_f =0$. Also, by lemma \ref{lema1}, if $e=f$ then $1_e 1_f = 1_e$ and hence $1_{e^{-1}}\delta_{e^{-1}}1_f\delta_f =  \alpha_{e^{-1}}(1_e)\delta_0=1_{e^{-1}}\delta_0$. Since $1_{e^{-1}}=1_r(e)$, (again by lemma \ref{lema1}), it follows that $1_{e^{-1}}\delta_{e^{-1}}1_e\delta_e=1_r(e) \delta_0$.

For the last relation, let $v\in E^0$ such that $0<\#\{e\in E^1:s(e)=v\}<\infty$. Note that in this case $X_v=\bigcup\limits_{e\in E^1:s(e)=v}^.X_e$. Then $$\sum\limits_{e\in E^1:s(e)=v} 1_e \delta_e 1_{e^{-1}}\delta_{e^{-1}} =\sum\limits_{e\in E^1:s(e)=v}1_e\delta_0=\left(\sum\limits_{e\in E^1:s(e)=v}1_e\right)\delta_0=1_v\delta_0.$$

So,  $\{1_e \delta_e, 1_{e^{-1}} \delta_{e^{-1}}: e\in E^1\}$ and $\{1_v\delta_0: v\in E^0\}$ satisfy conditions 1 to 4 of the definition of the Leavitt path algebra and hence, by the universal property, there exists a homomorphism $\varphi: L_K(E) \rightarrow D(X)\rtimes_\alpha\F$, such that, for all $e \in E^1$ and all $v\in E^0$, $\varphi(e)=1_e\delta_e$, $\varphi(e^*)=1_{e^{-1}}\delta_{e^{-1}}$ and $\varphi(v)=1_v\delta_0$ as desired.

\fim

Our next step is to show that the homomorphism $\varphi$ introduced above is bijective. To show that $\varphi$ is injective we will use the Graded Uniqueness Theorem of \cite{tomforde}, but in order to do so we first need to define a $\Z$-grading for $D(X)\rtimes_\alpha \F$. We do this below.

For each $p\in \F$, let $|p|:=m-n$, where $m$ is the number of generators (elements of $E^1$) of $p$ and $n$ is the number of inverses of generators of $p$. Define, for each $z\in \Z$, $A_z\subseteq D(X)\rtimes_\alpha\F$ as the $K$-linear span of $\{a_p\delta_p:a_p\in D_p \text{ and } |p|=z\}$, which is an additive subgroup of $D(X)\rtimes_\alpha \F$. It is straightforward to check that $A_zA_t\subseteq A_{z+t}$, and that $D(X)\rtimes_\alpha\F=\bigoplus\limits_{z\in \Z}A_z$ and so $\{A_z\}_{z\in \Z}$ is a $\Z$-grading of $D(X)\rtimes_\alpha\F$. 

\begin{theorem} The homomorphism $\varphi: L_K(E)\rightarrow D(X)\rtimes_\alpha \F$ introduced in the previous proposition is a $K$-isomorphism.
\end{theorem}

\demo First we show that $\varphi$ is injective. By the graded uniqueness theorem of \cite{tomforde} all we need to do is show that $\varphi$ is a graded ring homomorphism such that $\varphi(v) \neq 0$, for all $v\in E^0$. Recall, from \cite{tomforde}, that the Leavitt path algebra $L_K(E)$ has a $\Z$-grading $\{R_z\}_{\in \Z}$, where $R_z:=span\{ab^*:a,b\in W\cup E^0, |a|-|b|=z\}$, and we assume that $D(X)\rtimes_\alpha \F$ has a grading as defined above. Now, for each $ab^*\in R_z$, with $a,b\in W\cup\{E^0\}$, by the definition of $\varphi$, we have that $\varphi(ab^*)\in D_{ab^{-1}}\delta_{ab^{-1}}$, and since $|ab^{-1}|=|a|-|b|=z$ then $D_{ab^{-1}}\delta_{ab^{-1}}\subseteq A_z$ and hence $\varphi(R_z)\subseteq A_z$. Moreover, $\varphi(v)=1_v\delta_0\neq 0$ because $X_v\neq\emptyset$ for each $v\in E^0$. So, using the graded uniqueness theorem, we conclude that $\varphi$ is injective.

Before we prove that $\varphi$ is surjective we prove the following claim.

{\it Claim: For each $a\in W$, $\varphi(a)=1_a\delta_a$ and $\varphi(a^*)=1_{a^{-1}}\delta_{a^{-1}}$. Moreover, for each $a,b\in W$ with $a_{|a|}\neq b_{|b|}$ and $r(a)=r(b)$, $\varphi(ab^{*})=1_{ab^{-1}}\delta_{ab^{-1}}$. Hence, for each $p\in \F\setminus\{0\}$, $1_p\delta_p$ belongs to $Im(\varphi)$.}

Let $a=a_1...a_n\in W$. If $n=1$ then $\varphi(a)=1_a\delta_a$ by the definition of $\varphi$. Suppose $n\geq 2$ and $\varphi(a_2...a_n)=1_{a_2...a_n}\delta_{a_2...a_n}$. Then $$\varphi(a)=\varphi(a_1)\varphi(a_2...a_n)=1_{a_1}\delta_{a_1}1_{a_2...a_n}\delta_{a_2...a_n}=\alpha_{a_1}(\alpha_{a_1^{-1}}(1_{a_1})1_{a_2...a_n})\delta_{a_1a_2...a_n}=$$
$$=\alpha_{a_1}(1_{a_1^{-1}}1_{a_2...a_n})\delta_a=1_{a_1}1_{a_1a_2...a_n}\delta_a=1_a\delta_a,$$ where we used lemma \ref{lema2} to obtain the second to last equality.

In a similar way it follows that $\varphi(a^*)=1_{a^{-1}}\delta_{a^{-1}}$.

For $a,b\in W$ with $r(a)=r(b)$ and $a_{|a|}\neq b_{|b|}$, we have that $$\varphi(ab^{-1})=1_a\delta_a1_{b^{-1}}\delta_{b^{-1}}=\alpha_a(\alpha_{a^{-1}}(1_a)1_{b^{-1}})\delta_{ab^{-1}}=$$ $$=\alpha_a(1_{a^{-1}}1_{b^{-1}})\delta_{ab^{-1}}=1_a1_{ab^{-1}}\delta_{ab^{-1}}=1_{ab^{-1}}\delta_{ab^{-1}},$$ where we again made use of lemma \ref{lema2} to obtain the second to last equality.

This concludes the proof of the claim.

We can now prove that $\varphi$ is surjective, and for this it is enough to show that $D_p\delta_p \in Im(\varphi)$ for each $p\in \F$. 

First we show that $D_0\delta_0\in Im(\varphi)$. By linearity, it is enough to show that $1_v\delta_0\in Im(\varphi)$ and, for each $p\in \F\setminus\{0\}$, $1_p\delta_0\in Im(\varphi)$. We already know that $1_v\delta_0=\varphi(v)$. So, let $p\in \F\setminus\{0\}$. Notice that  $1_p\delta_0=1_p\delta_p1_{p^{-1}}\delta_{p^{-1}}$ and since $1_p\delta_p$ and $ 1_{p^{-1}}\delta_{p^{-1}} \in Im(\varphi)$, we conclude that $1_p\delta_0\in Im(\varphi)$.

To check that $D_p\delta_p\in Im(\varphi)$ it is enough to show, by linearity, that $1_p1_q\delta_q\in Im(\varphi)$ for each $p,q\in \F\setminus\{0\}$. But this is straightforward since $1_p1_q\delta_q=1_p\delta_01_q\delta_q$ and we already have that $1_p\delta_0\in Im(\varphi)$ and $1_q\delta_q\in Im(\varphi)$. 

We conclude that $\varphi$ is surjective as desired.

\fim

\section{Ideals in $D(X)\rtimes_\alpha \F$}

Now that we have realized the Leavitt path algebra of a graph as a partial skew group ring we wish to illustrate the power of the theory of partial skew group rings. We do so by providing new proofs of the Cuntz-Krieger uniqueness theorem and the simplicity criteria for Leavitt path algebras (see \cite{tomforde} for the original proofs), which relies only on the theory of partial skew group rings. Both proofs depend on the fundamental theorem \ref{TeorIdeal} below, which says that, for graphs that satisfy condition (L), that is, such that every closed path has an exit, any non-trivial ideal of $D(X)\rtimes_\alpha \F$ has non zero intersection with $D(X)$. To prove this theorem we need the following technical lemma:

\begin{lema}\label{MultiplyRightNonZero} Let $E$ be a graph that satisfies condition (L). If $b=b_1 b_2\ldots b_s \in W$ is a closed path and $x_b\in D_b$ is a non zero element, then there exists natural numbers $m,k\geq 1$, with $k\leq |b|$, and edges $t_1,...,t_k\in E^1$, such that $t_i\neq b_i$, for some $i$, and $x_b\cdot 1_{b^m t_1...t_k}\neq 0$.
\end{lema}

\demo 
Suppose first that there exists $N\in \N$ such that $x_b \cdot 1_{b^N} = x_b\cdot 1_{\underbrace{b...b}_{N\times}}=0$. Let $m$ be the greatest natural number that such that $x_b\cdot 1_{b^m}\neq 0$ (notice that $m\geq 1$, since $x_b=x_b\cdot 1_b$), so that $x_b\cdot 1_{b^{m+1}}=0$. Following the definition of the sets $X_{ab^{-1}}$, we have that  
$$ X_{b^m}=\bigcup\limits ^. X_{b^m t},$$ 
where the above disjoint union is taken over all $t\in W$ such that $s(t)=r(b)$ and $|t|=|b|$ or  $s(t)=r(b)$, $|t|<|b|$ and $r(t)$ is a sink.

From the description of $ X_{b^m}$ above, we conclude that there exists a $t=t_1...t_k\in W$ as above such that $x_b\cdot 1_{b^mt}\neq 0$. If $|t|<|b|$ then $r(t)$ is a sink, and in this case $t_k\neq b_k$, since $b_k$ is not a sink. If $|t|=|b|$ then $t\neq b$, since $x_b\cdot 1_{b^{m+1}}=0$, and hence $t_i\neq b_i$ for some $i$. 

Next suppose that $x_b\cdot 1_{b^N}\neq 0$ for all $N\in \N$. Since $x_b\in D_b$, we can write $x_b=\sum\limits_{j=1}^p\lambda_j1_{a_jc_j^{-1}}$, where $a_j, c_j\in W\cup \{0\}$. Choose $m\geq 1$ such that $m|b|\geq |a_j|$, for each $1\leq j\leq p$. Since $m|b|\geq |a_j|$, we infer that, for each $\xi \in X_b$, $(1_{b^m}1_{a_jc_j^{-1}})(\xi)$ depends only on $\xi_1,...,\xi_{m|b|}$, and since 
$$(x_b\cdot 1_{b^m})(\xi)=\sum\limits_{j=1}^p\lambda_j1_{b^m}(\xi)1_{a_jc_j^{-1}}(\xi)$$  
we also have that $(x_b\cdot 1_{b^m})(\xi)$ depends only on $\xi_1,...,\xi_{m|b|}$.

Now, since $x_b\cdot 1_{b^m}\neq 0$, there exists $\xi \in X_b$ such that $(x_b\cdot 1_{b^m})(\xi)\neq 0$ and, since the graph satisfies condition (L), there exists $t_i\in E^1$ with $s(t_i)=s(b_i)$ and $t_i\neq b_i$. Let $\mu \in X_b$ be such that $$\mu_1...\mu_{m|b|}\mu_{m|b|+i}=\xi_1...\xi_{m|b|}b_1...b_{i-1}t_i.$$ Then 
$$(x_b\cdot 1_{b^{m}b_1...b_{i-1}t_i})(\mu)=(x_b\cdot 1_{b^{m}})(\mu)=(x_b\cdot 1_{b^{m}})(\xi)\neq 0,$$ and so $x_b\cdot 1_{b^{m}b_1...b_{i-1}t_i}\neq 0$ as desired.

\fim

\begin{theorem}\label{TeorIdeal} Let $E$ be a graph which satisfies the condition (L). If $I$ is a non-zero ideal in $D(X)\rtimes_\alpha\F$ then $I\cap D(X) \neq 0$.
\end{theorem}

\demo Let $I$ be a non-zero ideal in $D(X)\rtimes_\alpha\F$ and $x$ a non-zero element of $I$. The proof we present below consists mostly of multiplying $x$ by appropriate elements of $D(X)\rtimes_\alpha\F$, so that eventually we obtain a new non-zero element in $I\cap D(X)$. In order to keep the proof organized we divide it in four steps.

{\it Step 1: There exists a non-zero element in $I$ of the form $\sum\limits_i x_{c_i}\delta_{c_i}$, where $c_i\in W\cup \{0\}$, and $c_i\neq c_j$ for each $i\neq j$}. 

 Let $x\in I$, $x\neq 0$. Then $x$ is a finite sum of the form $x=\sum\limits_i x_{a_ib_i^{-1}}\delta_{a_ib_i^{-1}}$, with $a_i,b_i\in W\cup\{0\}$, $a_ib_i^{-1}\neq a_jb_j^{-1}$ for $i\neq j$, and $x_{a_ib_i^{-1}}\neq 0$ for each $i$.  Choose $b_m$ among the elements $b_i$ appearing in the sum so that its length is equal to the greatest length of elements $b_i$ appearing in the sum, that is, so that $|b_m|=\displaystyle \max_i\{|b_i|\}$. Then $x\cdot 1_{b_m}\delta_{b_{m}}\neq 0$ because $x_{a_mb_m^{-1}}\delta_{a_mb_m^{-1}}1_{b_m}\delta_{b_m}=x_{a_mb_m^{-1}}\delta_{a_m}$, and $x_{a_mb_m^{-1}}\neq 0$. Also notice that $x\cdot1_{b_m}\delta_{b_m}=\sum\limits_i x_{c_i}\delta_{c_i}$ where $c_i\in W\cup\{0\}$, and so step 1 is proved.

{\it Step 2: There exists a non-zero element in $I$ of the form $x=x_0\delta_0+\sum\limits_{i=1}^n x_{d_i}\delta_{d_i}$, where $x_0\neq 0$, $d_i\in W$ for each $i$, $d_i\neq d_j$ if $i\neq j$, and if $i<j$ then $d_i$ is the beginning of $d_j$. Moreover $r(d_i)=s(d_i)=r(d_j)=s(d_j)$ for each $i,j$.}

Consider the element $x=\sum\limits_i x_{c_i}\delta_{c_i}\in I$ as in step 1, with $x_{c_i}\neq 0$ for each $i$ . Choose $c_n$ among the $c_i$ appearing in the sum so that $|c_n|=\displaystyle \max_i\{|c_i|\}$. If $c_n=0$ then $x=x_0\delta_0$, and the claim is proved. So, suppose $c_n\neq 0$. Note that $$y:=1_{c_n}\delta_0\cdot x=\sum\limits_i 1_{c_n}x_{c_i}\delta_{c_i}\neq 0,$$ since $1_{c_n}x_{c_n}\delta_{c_n}=x_{c_n}\delta_{c_n}$ and $x_{c_n}\neq 0$. Now,  for each $c_i$ with $|c_i|>0$,  $1_{c_n}x_{c_i}=1_{c_n}1_{c_i}x_{c_i}$ and since $1_{c_i}1_{c_n}=0$, unless $c_i$ is  the beginning of $c_n$, we have that $y$ is of the form $$y=y_0\delta_0+\left(\sum\limits_{i:c_i \text{ is the beginning of }c_n} y_{c_i}\delta_{c_i}\right),$$ where $y_{c_i}=x_{c_i}1_{c_n}$, for each $i$.

For one more reduction, we let $w=r(c_n)$ and multiply $y$ on the right by $1_w \delta_0$ to obtain
$$y\cdot 1_w\delta_0=y_01_w\delta_0+\left(\sum\limits_{i:c_i \text{ is the beginning of }c_n} \alpha_{c_i}(\alpha_{c_i^{-1}}(y_{c_i})1_w)\delta_{c_i}\right),$$ which, taking into consideration that $\alpha_{c_i^{-1}}(y_{c_i})1_w=\alpha_{c_i^{-1}}(y_{c_i})1_{c_i^{-1}}1_w=0$ if $r(c_i)\neq w$, can be written as $$y\cdot 1_w\delta_0=\left(\sum\limits_{c_i:r(c_i)=w, c_i \text{ is the beginning of } c_n}y_{c_i}\delta_{c_i}\right)+(y_01_w)\delta_0.$$

Next, we let $v=s(c_n)$, multiply $y\cdot 1_w\delta_0$ on the left by $1_v\delta_0$ and obtain
$$z:=1_v\delta_0\cdot y\cdot 1_w\delta_0=\left(\sum\limits_{c_i:r(c_i)=w, c_i \text{ is the beginning of } c_n}y_{c_i}\delta_{c_i}\right)+(y_01_v1_w)\delta_0.$$ Note that $z\neq 0$ since $y_{c_n}=x_{c_n}1_{c_n}=x_{c_n}\neq 0$.  

Now, if $y_01_v1_w\neq 0$, then $v=w$, i.e. $r(c_n)=s(c_n)$, and $z$ is the desired element. If $y_0 1_v 1_w=0$ then we ``pull" a non-zero element to the zero fiber. To do this, let $c$ be the element with the smallest length among the $c_i$´s corresponding to non-zero elements $y_{c_i}$ in the sum defining $z$ and consider the element $$1_{c^{-1}}\delta_{c^{-1}}\cdot z=\sum\limits_i\alpha_{c^{-1}}(\alpha_c(1_{c^{-1}})y_{c_i})\delta_{c^{-1}c_i}=\sum\limits_i\alpha_{c^{-1}}(y_{c_i})\delta_{c^{-1}c_i}.$$ Notice that $r(c^{-1}c_i)=r(c_i)=w=r(c)=s(c^{-1}c_i)$ and $1_{c^{-1}}\delta_{c^{-1}}\cdot z\neq 0$ since $y_c\neq 0$ and hence $\alpha_{c^{-1}}(y_c)\neq 0$. So $1_{c^{-1}}\delta_{c^{-1}}\cdot z$ is the desired element and this proves step 2.

{\it Step 3: Let $x=x_0\delta_0+\sum\limits_{i=1}^nx_{c_i}\delta_{c_i}$,  with $x_{c_n}\neq 0$, be a non-zero element of $I$ as described in step 2. Then there exists a non-zero $y\in I$ such that $y=y_0\delta_0$ or $y=y_0\delta_0+\sum\limits_{i=1}^k y_{c_i}\delta_{c_i}$ with $k<n$.}

First we rewrite $x$ as $\sum\limits_{i=0}^nx_{c_i}\delta_{c_i}$, where $c_0=0$, $x_{c_i}\neq 0$ for each $i$, and $r(c_i)=s(c_i)$ for each $c_i\neq 0$.
Next we choose $b$ among the $c_i's$ so that $|b|=\max\{|c_i|\}$. If $b=0$ then $x=x_0\delta_0$, and we are done. So, suppose $b\neq 0$ (so $b=c_n$) and write $b=b_1...b_{|b|}$. We will now construct a non-zero $y\in I$ as in the conclusion of step 3.

By lemma \ref{MultiplyRightNonZero}, there exists $m,k\geq 1$ and $t_1,...,t_k\in E^1$, with $t_i\neq b_i$ for some $i$, such that $(x_b\cdot1_{b^m t_1...t_k})\neq 0$. Consider the element $z$ of $I$ given by  
$$z:=(1_{b^{m}}\delta_0)\cdot x\cdot (1_{b^{m-1}t_1...t_k}\delta_0)=\sum\limits_{i=0}^nz_{c_i}\delta_{c_i},$$ where $z_{c_i}\delta_{c_i}=(1_{b^{m}}\delta_0)\cdot (x_{c_i}\delta_{c_i})\cdot 
(1_{b^{m-1}t_1...t_k}\delta_0)$. 

Notice that $z_{c_0}=0$, since $1_{b^{m}}1_{b^{m-1}t_1...t_k}=0$ and
$$z_{c_0}\delta_{c_0}=(1_{b^{m}}\delta_0)\cdot (x_0\delta_0)\cdot (1_{b^{m-1}t_1...t_k}\delta_0)=(x_{c_0}1_{b^{m}}1_{b^{m-1}t_1...t_k})\delta_0=0.$$ Hence 
$$z=\sum\limits_{i=1}^n (1_{b^{m}}\delta_0)\cdot (x_{c_i}\delta_{c_i})\cdot (1_{b^{m-1}t_1...t_k}\delta_0).$$
Furthermore, $z$ is a non-zero element, since 
$$z_{b}\delta_{b}=(1_{b^{m}}\delta_0)\cdot (x_{b}\delta_{b})\cdot (1_{b^{m-1}t_1...t_k}\delta_0)=
(x_b\cdot1_{\displaystyle b^{m}t_1...t_k})\delta_b\neq 0.$$

To obtain the desired element we have to ``pull" a non-zero element of $z$ to the zero fiber (as done in Step 2). To do so, let $c$ be the element with the smallest length among the $c_i$´s corresponding to non-zero elements $z_{c_i}$. Note that $1_{c^{-1}}\delta_{c^{-1}}\cdot z_c\delta_c=\alpha_{c^{-1}}(z_c)\delta_0$ and $\alpha_{c^{-1}}(z_c)\neq 0$ since $z_c\neq 0$, and hence $y:=1_{c^{-1}}\delta_{c^{-1}}\cdot z$ is the desired element. This proves step 3. 

{\it Step 4: Proof of the theorem.}

Let $x$ be a non-zero element of $I$ as in step 2, that is,  $x=x_0\delta_0+\sum\limits_{i=1}^nx_{c_i}\delta_{c_i}$. If $x_{c_i}=0$ for each $i$ then the theorem follows imediatily. If $x_{c_i}\neq 0$, for some $i$, then the theorem follows by finitely many applications of step 3.
\fim

We are now ready to prove the Cuntz-Krieger Uniqueness theorem, which will follow as a direct corollary of the result below.

\begin{proposition}\label{nonullintersection} Let $E$ be a graph that satisfies condition $(L)$ and let $I$ be a non-zero ideal in $D(X)\rtimes_\alpha \F$. Then there exists a vertex $v\in E^0$ such that $1_v\delta_0\in I$.
\end{proposition}

\demo 
Let $I$ be a non-zero ideal in $D(X)\rtimes_\alpha \F$. By the previous theorem there is a non-zero element $x_0\delta_0$ in $I$. Since $x_0\in D_0$, we can write $x_0$ as a linear combination of characteristic functions, namely  $x_0=\sum\limits_{i=1}^n\lambda_i1_{a_ib_i^{-1}}+\sum\limits_{j=1}^m\beta_j1_{v_j}$, where $a_i\in W$ and $b_i\in W\cup\{0\}$ (if $a_i=0$ then $1_{a_ib_i^{-1}}=1_{b_i}=1_{r(b_i)}$).

Let $v\in E^0$ be such that $1_vx_0\neq 0$. 

If $v$ is a sink, then 
$1_v1_{a_ib_i^{-1}} =0$ for each $i$, and then

$$0\neq 1_vx_0\delta_0=\sum\limits_{j=1}^n\beta_j1_v1_{v_j}\delta_0=\sum\limits_{j:v_j=v}\beta_j1_v\delta_0,$$ and so $1_v\delta_0\in I$.

% $a_i=0$ and $r(b_i) = v$, in which case $1_v1_{a_ib_i^{-1}}= 1_v$. We infer that 
%$$1_vx_0 \delta_0=\sum\limits_i\lambda_i1_v1_{a_ib_i^{-1}}\delta_0=\sum\limits_{i: a_i=0 \text{ and } r(b_i)=v}\lambda_i1_v1_{a_ib_i^{-1}}\delta_0=$$
%$$=\sum\limits_{i}\lambda_i1_v\delta_0=\left(\sum\limits_{i}\lambda_i\right)1_v\delta_0,$$
%and since $1_vx_0\neq 0 $, we conclude that $\sum\limits_{i: a_i=0 \text{ and } r(b_i)=v}\lambda_i\neq 0$ and hence $1_v\delta_0\in I$. 

Now we suppose that $v$ is not a sink. Let $m=max\{|a_i|:1\leq i\leq n\}$.

%If $m=0$ then $x_0=\sum\limits_{i=1}^n\lambda_i1_{b_i^{-1}}$. We can assume that one of the $b_i$´s is non-zero (call this element $b_0$), since if all the $b_i$´s are zero it is straightforward to find a vertex $v$ such that $1_v \delta_0 \in I$. Since $X_{b_i^{-1}} = X_{r(b_i)}$ we have that $$1_{r(b_0)}\delta_0 x_0 \delta_0 = \left( \sum\limits_{i:r(b_i)=r(b_0)}\lambda_i \right) 1_{r(b_0^{-1})}\delta_0$$ and hence $1_{r(b_0^{-1})} \delta_0 \in I$.

%Next we deal with the case when $m>0$. 
Recall that we can write $X_v$ as
$$X_v=\bigcup\limits^.X_c,$$ where the above disjoint union is taken over all $c\in W$ such that $s(c) = v$ and $|c|=m$ or $s(c)=v$, $|c|<m$ and $r(c)$ is a sink. From this we deduce that there is a $c$ in the above union index set such that $1_cx_0\neq 0$. Now, if $1_c1_{a_ib_i^{-1}}\neq 0$ then $a_i$ is the beginning of $c$, and then $1_c1_{a_ib_i^{-1}}=1_c1_{a_i}=1_c$. Moreover, if $1_c1_{v_j}\neq 0$ then $1_c1_{v_j}=1_c$.

%, or $a_i=0$ and then $1_c1_{a_ib_i^{-1}}=1_c1_{b_i^{-1}}=1_c$. 

With this in mind, we obtain that 
$$0\neq 1_cx_0 \delta_0 =\sum\limits_{i=1}^n\lambda_i1_c1_{a_ib_i^{-1}}\delta_0+\sum\limits_{j=1}^m\beta_j1_c1_{v_j}\delta_0=$$
$$=\sum\limits_{i:1_c1_{a_ib_i^{-1}}\neq 0}\lambda_i1_c1_{a_ib_i^{-1}}\delta_0+\sum\limits_{j:1_c1_{v_j}\neq 0}\beta_j1_c1_{v_j}\delta_0=$$
$$=\sum\limits_{i:1_c1_{a_ib_i^{-1}}\neq 0}\lambda_i1_c\delta_0+\sum\limits_{j:1_c1_{v_j}\neq 0}\beta_j1_c\delta_0=\left(\sum\limits_{i:1_c1_{a_ib_i^{-1}}\neq 0}\lambda_i+\sum\limits_{j:1_c1_{v_j}\beta_j\neq 0}\beta_j\right)1_c\delta_0,$$ and hence $1_c\delta_0$ is a non-zero element in $I$. Since $I$ is an ideal, we have that  $1_{c^{-1}}\delta_0=1_{c^{-1}}\delta_{c^{-1}}\cdot1_c\delta_0\cdot1_c\delta_c$ also belongs to $I$, and since $1_{c^{-1}}=1_{r(c)}$ this proves the proposition.
\fim

\begin{corolario}(Cuntz-Krieger Uniqueness theorem) Let $E$ be a graph that satisfies condition $(L)$. If $\phi:D(X)\rtimes_\alpha\F\rightarrow B$ is a $K$-homomorphism such that $\phi(1_v\delta_0)\neq 0$ for each $v\in E^0$, then $\phi$ is injective.
\end{corolario}

\begin{obs} The above theorem was first proved (via a different approach) in \cite{tomforde}.
\end{obs}

Our final goal is to present a new proof for the simplicity criteria for Leavitt path algebras. For the reader convinience we recall the definition of hereditary and saturated subsets of $E^0$ below.

\begin{definition}
\begin{enumerate}
\item A subset $H\subseteq E^0$ is said to be hereditary if for each $v\in H$ such that $v=s(e)$, for some $e\in E^1$, we have that $r(e)\in H$. 
\item A subset $H\subseteq E^0$ is called saturated if whenever $v\in E^0$ is such that $0<|s^{-1}(v)|<\infty$, then $\{r(e):e\in E^1 \text{ and }s(e)=v\}\subseteq H$ implies $v\in H$.
\end{enumerate}
\end{definition}

Let I be an ideal in $D(X)\rtimes_\alpha\F$ and define $H_I=\{v\in E^0: 1_v\delta_0 \in I\}$. By proposition \ref{nonullintersection}, if $I\neq 0$ and $E$ satisfies condition $(L)$ then $H_I$ is nonempty. The following result appeared first in \cite{arandapino1} (lemma 2.3) in the Leavitt path algebras context. Below we give the partial skew group ring version of it.

\begin{proposition} Let I be and ideal in $D(X)\rtimes_\alpha\F$. Then the set $H_I$ is saturated and hereditary.
\end{proposition}

\demo First we prove that $H_I$ is hereditary. For this, suppose that $v=s(e)\in H_I$. Notice that since $1_{s(e)}\delta_0\in I$ we have that $1_e\delta_e=1_{s(e)}\delta_0\cdot 1_e\delta_e\in I$ and so $1_{e^{-1}}\delta_0=1_{e^{-1}}\delta_{e^{-1}}\cdot1_e\delta_e\in I$. We conclude that $1_{r(e)}\delta_0\in I$ (since $1_{r(e)}=1_{e^{-1}}$) and hence $r(e)\in H_I$ as desired.

Now we prove that $H_I$ is saturated. Let $v\in E^0$ be such that $0<|s^{-1}(v)|<\infty$ and suppose $\{r(e):e\in E^1 \text{ and }s(e)=v\}\subseteq H_I$. Let $e\in E^1$ be such that $s(e)=v$. Since $r(e)\in H_I$ we have that $1_e\delta_e=1_e\delta_e\cdot 1_{r(e)}\delta_0\in I$, and so $1_e\delta_0=1_e\delta_e\cdot 1_{e^{-1}}\delta_{e^{-1}}\in I$. But his implies that $1_v\delta_0=\sum\limits_{e\in E^1:s(e)=v}1_e\delta_0\in I$, since $1_v=\sum\limits_{e\in E^1:s(e)=v}1_e$, and hence $v\in H_I$ as desired. \fim 

We can now prove the simplicity criteria for Leavitt path algebras (originally proved in \cite{arandapino1, tomforde} for $L_K(E)$) using only partial skew group ring theory.

\begin{proposition} Let $E$ be a graph that satisfies condition $(L)$, and suppose that the only saturated and hereditary subsets of $E^0$ are $\emptyset$ and $E^0$. Then $D(X)\rtimes_\alpha\F$ is simple.
\end{proposition}

\demo Let $I$ be a non-zero ideal in $D(X)\rtimes_\alpha\F$, and let $H_I$ as above. By the previous proposition, $H_I$ is saturated and hereditary. By proposition \ref{nonullintersection}, $H_I$ is nonempty, and hence $H_I=E^0$. We conclude that $1_v\delta_0\in I$ for each $v\in E^0$. 

Now, for each $p\in \F\setminus\{0\}$, there is a $v\in E^0$ such that $1_v1_p=1_p$ (recall that if $a\in W$ then $1_{s(a)}1_{ab^{-1}}=1_{ab^{-1}}$, and $1_{r(a)}1_{a^{-1}}=1_{a^{-1}}$). So, let $p\in \F\setminus\{0\}$, $q\in\F$, and let $v\in E^0$ be such that $1_v1_p=1_p$. Then  $1_p1_q\delta_q=1_v\delta_0\cdot1_p1_q\delta_q\in I$ and, since $I$ is an ideal, it follows that $D_p\delta_p\subseteq I$ for all $p\in \F$, and hence $D(X)\rtimes_\alpha\F\subseteq I$.
\fim

\vspace{1.5pc}

D. Goncalves, Departamento de Matemática, Universidade Federal de Santa Catarina, Florianópolis, 88040-900, Brasil

Email: daemig@gmail.com

\vspace{0.5pc}
D. Royer, Departamento de Matemática, Universidade Federal de Santa Catarina, Florianópolis, 88040-900, Brasil

Email: royer@mtm.ufsc.br
\vspace{0.5pc}

\end{document}